\documentclass[letterpaper, 10 pt, conference]{ieeeconf}  

\IEEEoverridecommandlockouts                              
\usepackage[utf8]{inputenc}
\usepackage[T1]{fontenc}
\let\labelindent\relax
\usepackage{enumitem}
\usepackage{amsmath, amssymb}
\usepackage{semantic}
\usepackage{nicematrix}
\NiceMatrixOptions{letter-for-dotted-lines=V}
\usepackage{arydshln}
\usepackage{algpseudocode}
\usepackage{algorithm}
\usepackage{graphicx}
\usepackage{comment}
\usepackage{balance}
\usepackage{optidef}
\usepackage{mathrsfs}
\usepackage{varwidth}
\newtheorem{theorem}{Theorem}
\newtheorem{proposition}{Proposition}
\newtheorem{lemma}{Lemma}
\newtheorem{definition}{Definition}
\newtheorem{remark}{Remark}
\newtheorem{assumption}{Assumption}
\newtheorem{problem}{Problem}

\newcommand{\R}{\mathbb{R}}
\newcommand{\Rn}{\mathbb{R}^n}
\newcommand{\Rnn}{\mathbb{R}^{n \times n}}
\newcommand{\qed}{\hfill $\blacksquare$}

\newcommand{\bac}[1]{\textcolor{blue}{#1}}

\title{\LARGE \bf
Distributed Optimization of Average Consensus Containment with Multiple Stationary Leaders
}

\author{Sushobhan Chatterjee and Rachel Kalpana Kalaimani
\thanks{*This work has been partially supported by DST-INSPIRE Faculty Grant, Department of Science and Technology (DST), Govt. of India (ELE/16-17/333/DSTX/RACH)}
\thanks{The authors are with the Department of Electrical Engineering, Indian Institute of Technology Madras, Chennai, India - 600036. \texttt{(email: chatterjeesushobhan8@gmail.com; rachel@ee.iitm.ac.in)}.}%
}

\begin{document}

\maketitle
\thispagestyle{empty}
\pagestyle{empty}

\begin{abstract}

In this paper, we consider the problem of containment control of multi-agent systems with multiple stationary leaders, interacting over a directed network. While, containment control refers to just ensuring that the follower agents reach the convex hull of the leaders' states, we focus on the problem where the followers achieve a consensus to the average values of the leaders' states. We propose an algorithm 
that can be implemented in a distributed manner to achieve the above consensus among followers. Next we optimize the convergence rate of the followers to the average consensus by proper choice of weights for the interaction graph.  This optimization is also performed in a distributed manner using Alternating Direction Method of Multipliers (ADMM). Finally, we complement our results by illustrating them with numerical examples.\\  
\end{abstract}


\section{INTRODUCTION}
The idea of achieving a common objective in a system with multiple agents using distributed and co-operative control  has garnered a lot of interest among research community for a long time.  A fundamental aspect of the distributed co-operative control is that the agents achieve the collective objective via local interaction among neighbours, thereby making the framework robust to communication failures, highly adaptable and cost effective.

A well addressed problem in this multi-agent system framework is the consensus problem, where the agents are expected to agree on a common value \cite{ren2005consensus}, \cite{olfati2007consensus}. 
 Consensus with a  single leader is addressed in \cite{jadbabaie2003coordination}. 
%
When there are multiple leaders the problem is referred to as containment control where the objective is to bring the remaining follower agents from any arbitrary state, within the convex hull spanned by the leaders' state \cite{xiao2006consensus}. This problem is also motivated from nature such as gravel ants forming a boundary over caterpillars and transporting along with them due to their sugary extract. Practical applications include shaping of social dynamics using a set of guiding agents, military applications, safe disposal of hazardous waste etc. \cite{cao2010distributed}. Extensive research is done in this area such as analysing dynamic leaders using PDE in \cite{ji2008containment}, assuming directed static topology \cite{liu2012necessary}, time-varying topologies \cite{cao2009containment}, agents with higher order dynamics and 
impact of time delays \cite{li2018distributed}. 

While the above containment algorithms only focus on general containment of followers, i.e. each follower achieving an arbitrary state, different from each other, within the convex hull of the leaders, the primary objective of this paper is to achieve a consensus among the followers regarding the final destination within the convex hull of the leaders.
The idea of achieving consensus in a multiple leader containment scenario is motivated by various applications. One is opinion consensus of a group, driven by a couple of external agents/factors, in a direction that gets aligned with their own priors. Another application is a cooperative attack on a static or, moving target, where the target is encircled by UAVs (leaders) such that the target's live location coincides with the average of coordinates of UAVs. Then the missiles (followers) can easily home in to the target without the need for any individual radar tracking.      

Motivated by the above applications for consensus containment, we provide an algorithm where the followers reach a consensus to the average of the leaders' state. \cite{khan2010higher} investigates consensus of followers inside the convex hull of stationary leaders, where a distributed control law has been designed by formulating a centralized multi-objective optimization (MOP) with a trade-off between performance error and convergence rate. 
This results in a significant asymptotic drift in the final consensus value due to non-convergence to the exact solution on account of trade-offs.
We propose a distributed control law which ensures a fast convergence of the followers to the consensus within the convex hull spanned by the leaders over a directed network, without any aforementioned asymptotic drift. The fast convergence is designed using a distributed optimization algorithm unlike the centralized approach in \cite{khan2010higher}. 
Distributed optimization of convergence rate to consensus without leaders in a multi-agent, undirected framework is addressed in \cite{rokade2020distributed}. 

We summarize our contribution as follows.
\begin{enumerate}
    \item We consider a multi-agent system with multiple stationary leaders and followers communicating with each other over a directed network topology. We propose a distributed algorithm which ensures consensus of followers to the average of leaders' states. 
    \item  We formulate an optimization problem to maximize the rate of convergence of follower agents to the consensus value. We propose an algorithm based on Alternating Direction Method of Multipliers (ADMM) that solves the above optimization problem in a distributed manner.
    \item We illustrate the results using numerical examples. We compare the results of the distributed optimization to the centralized optimization. 
\end{enumerate}

The outline of the paper is laid out as follows. Following the introduction, Section II focuses on basic mathematical preliminaries. In Section III, the two problems are formally discussed. Algorithms for solving these problems are proposed in Section IV. Section V provides a numerical example to elucidate the validity of the proposed algorithms. Finally, we give concluding remarks in Section VI. The proofs of some of the results are in Appendix.

\section{PRELIMINARIES}
In this section, we briefly introduce some preliminaries pertaining to directed graphs and convex analysis.
\subsection{Notations}
Let $\Rn$ be the $n$-dimensional real vector and $\Rnn$ be the real square matrix of order $n$. Let $I_n$ be the $n \times n$ identity matrix, $0_n$ be the $n \times n$ matrix with all entries as 0 and the vector $1_n = [1,\hdots, 1]^T$ with dimension $n$. For any vector $x$ : $x_i$ or $\lbrack x \rbrack_i$ denotes its $i^{th}$ element, $\overline{x}$ denotes its average value, and $\Vert x \rVert_2$ denotes its 2-norm. For a set S, |S| denotes its cardinality. For a matrix $P$ : $\Lambda(P)$ and $\rho(P)$ denote its eigenvalue set and spectral radius, respectively, $P^{ij}$ or $(P)^{ij}$ denotes its $(i, j)^{th}$ element and $\lVert P \rVert_F$ denotes its Frobenius-norm.  
 \subsection{Graph Theory}
Let $\mathcal{G} := (\mathcal{V},\mathcal{E},A)$ denote a directed graph, where $\mathcal{V} = [v_1,v_2,...,v_n]$ represents the vertex set, $\mathcal{E} \subseteq \mathcal{V}\times \mathcal{V}$ denotes the directed edge set and $A$ = [$A^{ij}$] $\in \Rnn$ denotes the weighted adjacency matrix defined as follows.
A directed edge is an ordered pair of distinct vertices $(v_j,v_i)$ such that, $i^{th}$ node (child) can access state information of $j^{th}$ node (parent). The weight associated with the edge $(v_j,v_i)$ is the entry $A^{ij} \geq 0$, in $A$.  Also we assume $A^{ii} > 0\ \forall\ i = \{1,2,...,n\}$.\\ $L$ = [$L^{ij}$] $\in \Rnn$ denotes the laplacian matrix defined as: $L^{ij} = \overset{k=n}{\underset{k=1, k\neq i}{\Sigma}}A^{ik}\ \forall\ i=j$, and $L^{ij}=-A^{ij}\ \forall\ i\neq j$. For a vertex $v_i$, the set of in-neighbours is defined by $\mathcal{N}^{in}_i=\big\{v_j\ |\ (v_j,v_i)\in \mathcal{E}\big\}$ and out-neighbours by $\mathcal{N}^{out}_i=\big\{v_j\ |\ (v_i,v_j)\in \mathcal{E}\big\}$.
For a graph $\mathcal{G} := (\mathcal{V},\mathcal{E},A)$, the subgraph $\mathcal{G}^1$ is the graph induced by a vertex set $\mathcal{V}^1 \subseteq \mathcal{V}$ and edge set $\mathcal{E}^1=\big\{(v_i,v_j)\in \mathcal{E}\ \forall\ v_i,v_j\in\mathcal{V}^1\big\}.$\\
A directed graph $\mathcal{G}$ is said to be strongly connected if there exists a sequence of directed paths between any two distinct pair of vertices $\{v_j,v_i\} \in \mathcal{V}$, starting at $v_j$ and ending at $v_i$. 


 \subsection{Convex Analysis}   
 A set $D \subset \mathbb{R}^q$ is said to be convex if $\exists\ x,\ y \in D\ s.t.$\\ $(1-\theta)x + \theta y \in D\ \forall\ \theta \in [0,1]$. 
 \begin{definition}
\cite{rockafellar1997convex}
 The convex hull of a finite set of points $\{p_1,...,p_n\} \in \mathbb{R}^q$ is the smallest convex set containing all points $p_i$, $i = 1,2,\hdots,n$ denoted by \textbf{conv}$\{p_1,...,p_n\} = \bigg\{\overset{n}{\underset{i=1}{\Sigma}} \theta_ip_i\ |\ \theta_i \in \mathbb{R},\ \theta_i \geq 0,\ \overset{n}{\underset{i=1}{\Sigma}} \theta_i = 1\bigg\}$. 
\end{definition}

\section{PROBLEM FORMULATION} 

We first explain the framework of a multi-agent system with more than one leader. 
Consider a group of $n$ agents communicating with each other.
The communication pattern is depicted by a directed graph $\mathcal{G}$. The vertices of the graph correspond to the agents and 
the state of an agent $j$ is accessible by agent $i$ only if there is a directed edge $(v_j,v_i)$ in the graph. 
An agent is designated to be a leader if its in-neighbour set is empty. The rest of the agents are termed as followers. 
Let $\mathcal{L}$ and $\mathcal{F}$ denote respectively the index set of leaders and followers.
Let $\mathcal{G^F}$ be the subgraph induced by the vertex set of the followers, $v_i\ \forall\ i \in \mathcal{F}$ and $\mathcal{G^{FL}}$ is a subgraph with vertex set $\mathcal{V}$ and edge set given by $\mathcal{E^{FL}}=\big\{(v_j,v_i)\in \mathcal{E}\ \forall\ i\in\mathcal{F},\ j\in\mathcal{L}\big\}$.
Let $x_i\in\R$ denote the state of each agent and assume that they follow the following discrete time dynamics.
\begin{equation} \label{eqn:system}
x_i(k+1)=\alpha_ix_i(k) + u_i(k),\end{equation}
where $u_i(k)$ is the control input to agent $i$ at $k^{th}$ iteration.  The following are our assumptions on the agents.
\begin{assumption} \label{assume}
\begin{enumerate}[label=(\alph*)]
    \item The leaders are stationary, i.e., $\forall\ k \geq 0\ \&\ i \in \mathcal{L},\ u_i(k) = 0\ \&\ \alpha_i = 1$, and hence $x_i(k+1)=x_i(k)$.
    \item Every leader has at least one follower connected to it.
    \item $\mathcal{G^F}$ is strongly connected. 
\end{enumerate}
\end{assumption}
\subsection{Average Consensus Containment} 
Containment control refers to the problem of designing an input for the followers such that followers reach the convex hull spanned by the stationary leaders \cite{cao2009containment}, \cite{liu2012necessary}. 
Containment control does not guarantee a consensus within the convex hull. We focus on the problem where followers reach a consensus to average of leaders' states. We term this problem as \textit{average consensus containment control}. We explain later how the existing containment protocol can be modelled to achieve average consensus containment using a centralized approach. 
We propose to address this problem using a distributed approach.
Our problem is formulated below.
\begin{problem}
Consider a set of $n$ agents with $m$ leaders and $n-m$ followers interacting over a directed graph satisfying \textit{Assumption 1} with dynamics given in \eqref{eqn:system}. 
Find a control input that can be designed and implemented in a distributed manner,  such that the follower agents state converges to the average of the stationary leaders state.
\begin{equation}\label{eqn:convergence}
\underset{k \rightarrow \infty}{\lim}x_i(k)\rightarrow \frac{1}{m} {\underset{j\in\mathcal{L}}{\sum}}x_j(0)\ \forall\ i \in \mathcal{F} 
\end{equation}
\end{problem}

\subsection{Optimizing the convergence rate} 
Our next objective is to maximize the rate of convergence of the followers to the consensus mentioned in Problem 1. This can be achieved by suitably choosing the edge weights of the network over which the agents are communicating. Based on the results that we obtain for Problem 1, we first formulate the fastest convergence rate problem as an optimization problem. Then we attempt to solve that in a distributed manner. This is formulated as follows.
\begin{problem}
Consider Problem 1 with the mentioned assumptions and network model. We solve the following:
\begin{enumerate}
    \item Formulate an optimization problem that maximizes the rate of convergence given in \eqref{eqn:convergence}.
    \item Propose to solve this optimization problem in a distributed manner.
\end{enumerate}

\end{problem}


\section{MAIN RESULTS}
\subsection{Average Consensus Containment}
In this section, we first briefly explain the assumptions and protocol that ensures containment control. Then we introduce our algorithm which achieves average consensus containment that can be implemented in a distributed manner. 


Consider a multi-agent network of $n$ agents with $m$ stationary leaders and $n-m$ followers, having dynamics as given in \eqref{eqn:system} interacting over a network depicted by a graph $\mathcal{G}$. Assume that the agents use the following update protocol,
\begin{equation}
u_i(k) = 
     \begin{cases}
      \underset{j \in \mathcal{N}^{in}_i}{\Sigma}\alpha_{ij}x_j(k), &\quad i \in \mathcal{F}\\
     \quad\quad\quad 0, &\quad i \in \mathcal{L}
     \end{cases} \null \label{eqn:GC_u}
\end{equation}
for some suitable choice of weights $\alpha_{ij}$
Let $x_L(k) \in \mathbb{R}^m$ and $x_F(k) \in \mathbb{R}^{n-m}$ refer to the vector of all leader and follower states respectively at $k^{th}$ iteration. Based on the weights in \eqref{eqn:GC_u} and the system dynamics in \eqref{eqn:system}, we define a matrix $A\in\Rnn$, where $A^{ij}=\alpha_{ij}$, for $i\neq j$ and $A^{ii}=\alpha_i$ which represents a weighted adjacency matrix corresponding to graph $\mathcal{G}$.  
The agent dynamics can be written in matrix form as follows. 
\begin{equation} \label{eqn:Apartition}
\begin{bmatrix}x_F(k+1)\\x_L(k+1)\end{bmatrix} = \begin{bmatrix}
A_1 & A_2\\
0 & I_m
\end{bmatrix}\begin{bmatrix}x_F(k)\\x_L(k)\end{bmatrix}, A=\begin{bmatrix}
A_1 & A_2\\
0 & I_m
\end{bmatrix}
\end{equation}
The following result states the conditions for containment of followers within the convex hull of leader's states for a suitable choice of $A_1$ and $A_2$. 
\begin{proposition}
\textit{(\cite{wang2019necessary}, Theorem 9) Consider a set of multi-agents  satisfying Assumption \ref{assume} following the protocol in \eqref{eqn:Apartition}.  
Assume $A_1=I_{n-m}-\alpha L_1$, $A_2=-\alpha L_2$, where $L_1$ and $L_2$ are the weighted laplacian matrices associated with subgraphs $\mathcal{G^F}$ and $\mathcal{G^{FL}}$ respectively. The followers will converge to the convex hull formed by stationary leaders if and only if step-size satisfies; $\alpha < \underset{\lambda_i\in \Lambda(L_1)}{min}\ \frac{2Re(\lambda_i)}{Re^2(\lambda_i) + Im^2(\lambda_i)}$, and the final state of all the followers are given by $-L_1^{-1}L_2x_L(0)$}. \qed
\end{proposition}
From the above proposition, we observe that the follower states converge to some value within the convex hull formed the stationary leader's states and do not achieve any consensus.
We next proceed to design an update protocol for the \textit{average consensus containment problem}.
We propose Algorithm \ref{algo:avgcont} for average consensus containment under suitable assumptions on the matrix $A$. The existing Push sum algorithm \cite{kempe2003gossip} is adapted to the current problem. The convergence of the algorithm is proved in Theorem \ref{theo:algoconvg}.
\begin{remark} ({\em Stopping Criterion for Algorithm \ref{algo:avgcont}}) Each of the agent $i \in \mathcal{F}$ fixes a small, arbitrary tolerance value, $\gamma > 0$, and calculates its \textit{successive iterate error} $\{e_i(k)\}$, at each iteration $k$ as
\begin{equation}
    e_i(k) = x_i(k) - x_i(k-1) 
\end{equation}
which indicates the extent of asymptotic convergence of the agents trajectory. The stopping criterion for Algorithm \ref{algo:avgcont} is when
\begin{equation}
    E_i(k) = \big|e_i(k)\big| \leq \gamma\ \forall\ i \in \mathcal{F}
\end{equation}
\end{remark}

\begin{algorithm}
  \caption{Distributed Average Consensus Containment}
  \label{algo:avgcont}
  \begin{algorithmic}
    \State \textbf{Given :} $A$, Initial condition $x(0)$ of all agents.
    \State \textbf{Assumptions :} $A_1$ and $A_2$ defined in \eqref{eqn:Apartition} are column stochastic and agents satisfy Assumption \ref{assume}.
    \State \textbf{Initialize :} $s(0) = x(0),\ w(0) = 1_n$
    \State \textbf{Iterate :}
    \For{ $k \geq 0$ } 
    \State \textbf{Step 1 :} \textbf{Exchange Values :}  $\forall\ i \in \mathcal{F}$,
    \State   Receive $s_j(k),\ w_j(k)$ from all in-neighbours $j \in \mathcal{N}^{in}_i$  
    \State \textbf{Step 2 :} \textbf{Update sum and weight vectors :} 
    \State $s(k+1) = As(k)$
    \State $w(k+1) = Aw(k)$
    \State \textbf{Step 3 :} \textbf{Update state vector :} $\forall\ i \in \mathcal{F}$
    \State $x_i(k+1) = \frac{s_i(k+1)}{w_i(k+1)}$ 
    \EndFor  
  \end{algorithmic}
\end{algorithm}


\begin{remark} ({\em Algorithm 1 implemented in a distributed and parallel setting}) 
Each agent $i$ maintains three variables, $s_i(k)$, $w_i(k)$ and $x_i(k)$ at each iteration $k$. At every iteration, each agent updates its variables $s_i$ and $w_i$ by exchanging data only with its neighbours as indicated in Step 1 and 2. Furthermore, every agent updates the variable $x_i$ as indicated in Step 3, parallely. 
\end{remark}
\begin{theorem}\label{theo:algoconvg}
\textit{The average consensus containment protocol proposed in Algorithm 1 ensures that followers converge to a consensus which is the average of the leaders' state values.} 
\end{theorem} 
The proof is given in Appendix.\qed 
\subsection{Optimizing the convergence rate}     In this section we address Problem 2, where we maximize the rate of convergence of the follower agents to consensus which is the average of stationary leaders' states. We first identify the parameter that influences the rate of convergence and then formulate the relevant optimization problem. 
The following lemma characterises the rate of convergence of the follower agents in terms of the weighted adjacency matrix.

\begin{lemma}
Consider a multi-agent system following dynamics given in \eqref{eqn:system} and satisfying Assumption \ref{assume}. Assume that the agents are implementing Algorithm \ref{algo:avgcont}.  The rate of convergence is characterised by $\rho\big(A_1 - \frac{1}{n-m}\mathbf{11}^T\big)$. \qed
\end{lemma}
The proof of the lemma is given in Appendix.\\ 
From the above lemma, it is clear that we have to optimize the weights of the follower interaction, i.e., entries of $A_1$, so as to maximize the convergence rate. For this, we formulate the following optimization problem with appropriate constraints. 
\begin{mini*}|1|
 {\substack{A_1}}{\rho\bigg(A_1 - \frac{1}{n-m}\mathbf{11}^T\bigg)}{}{}
 \label{eqn:P1}
 \addConstraint{A_1\mathbf{1} = \mathbf{1},\ A_1^T\mathbf{1} = \mathbf{1}} \tag{$A$}
 \addConstraint{(A_1)^{ij} = 0,\ (i,j) \notin \mathcal{E^F}}
\end{mini*}
where $\mathcal{E^F}$ is the edge set of the subgraph $\mathcal{G^F}$ and $\mathbf{1} = 1_{n-m}$. Here, $A_1$ should be column stochastic because of the assumption in Algorithm \ref{algo:avgcont}, and 
the requirement of row stochasticity has been explained later. The third equality is the topological constraint imposed by the communication pattern of the follower agents. 

The above minimization problem is not convex due to $\rho\big(A_1 - \frac{1}{n-m} \mathbf{11}^T\big)$ being a non-convex function of $A_1$  \cite{overton1988minimizing}.
The problem is modified by using a convex function $\big\lVert A_1 - \frac{1}{n-m}\mathbf{11}^T \big\rVert_2$ in place of the objective function $\rho\big(A_1 - \frac{1}{n-m}\mathbf{11}^T\big)$. The relaxed optimization problem is as follows.
\begin{mini*}|1|
 {\substack{A_1}}{\bigg\lVert A_1 - \frac{1}{n-m}\mathbf{11}^T \bigg\rVert_2}{}{}
 \label{eqn:P2}
 \addConstraint{A_1\mathbf{1} = \mathbf{1},\ A_1^T\mathbf{1} = \mathbf{1}} \tag{$B$}
 \addConstraint{(A_1)^{ij} = 0,\ (i,j) \notin \mathcal{E^F}}
\end{mini*}
Since, for a directed graph, $A_1 \neq A_1^T$, $\rho\big(A_1 - \frac{1}{n-m}\mathbf{11}^T\big) \leq \lVert A_1 - \frac{1}{n-m}\mathbf{11}^T \rVert_2$.  
Suppose $A_1^{*}$ is the solution to above optimization problem \eqref{eqn:P2}, then with $A_1^{*}$ being doubly stochastic it can be shown that $\lVert A_1^{*} - \frac{1}{n-m}\mathbf{11}^T \rVert_2 < 1$ (\cite{rokade2020distributed}, Lemma 3), and hence $\rho\big(A_1^{*} - \frac{1}{n-m}\mathbf{11}^T\big) < 1$.

Note that the above problem can be solved as a centralized optimization problem using any of the standard optimization solvers. Since our proposed Algorithm \ref{algo:avgcont} is implemented in a distributed manner, we propose to solve the above problem \eqref{eqn:P2} also in a distributed manner using ADMM over directed graphs. For this purpose, we modify the problem as follows.
\begin{mini*}|1|
 {\substack{(A_1)_i\ \forall\ i \in \mathcal{F}}}{\overset{n-m}{\underset{i=1}{\Sigma}} \frac{\big\lVert (A_1)_i-\frac{1}{n-m}\mathbf{11}^T \big\rVert_2}{n-m}}{}{}
 \label{eqn:P3}
 \addConstraint{(A_1)_i = (A_1)_j,\ (i,j) \in \mathcal{E^F}}
 \addConstraint{(A_1)_i\mathbf{1} = \mathbf{1},\ (A_1)_i^T\mathbf{1} = \mathbf{1}} \tag{$C$}
 \addConstraint{(A_1)_i^{ij} = 0,\ (i,j) \notin \mathcal{E^F}}
\end{mini*}
In problem \eqref{eqn:P3}, 
each of the matrices $(A_1)_i$ are copies of the optimization variable $A_1$ maintained by each agent $i \in \mathcal{F}$, and are their estimates of the centralized optimal solution $A_1^{*}$. Now, since $\mathcal{G^F}$  is strongly connected, $(A_1)_i = (A_1)_j,\ (i,j) \in \mathcal{E^F}$ implies that all $(A_1)_i$'s are equal. Hence, we conclude that \eqref{eqn:P2} and \eqref{eqn:P3} are equivalent. 
Since we want the constraints to be decoupled for a distributed approach, 
we ensure that each $(A_1)_i$ has the graph topological constraints only on its $i^{th}$ row and not the entire matrix. This is ensured in the last constraint of the above problem.\\
Now \eqref{eqn:P3} can be implemented in a distributed manner, but parallel operation across all the agents is not possible due to the fact that $(A_1)_i$'s of different agents all come up in the same constraint, $(A_1)_i = (A_1)_j,\ (i,j) \in \mathcal{E^F}$ in \eqref{eqn:P3}. In order to circumvent this issue, we formulate it using Fenchel duality (\cite{notarstefano2019distributed}, Section 3.1.2). We introduce an auxiliary primal variable $Z \in \mathbb{R}^{(n-m) \times (n-m)}$ as follows, which decouples the consensus constraint by enforcing cohesion among all copies.
\begin{mini*}|1|
 {\substack{(A_1)_i\ \forall\ i \in \mathcal{F}}}{\overset{n-m}{\underset{i=1}{\Sigma}} \frac{\big\lVert (A_1)_i-\frac{1}{n-m}\mathbf{11}^T \big\rVert_2}{n-m}}{}{}
 \label{eqn:P4}
 \addConstraint{(A_1)_i = Z,\ i \in \mathcal{F}}
 \addConstraint{(A_1)_i\mathbf{1} = \mathbf{1},\ (A_1)_i^T\mathbf{1} = \mathbf{1}} \tag{$D$}
 \addConstraint{(A_1)_i^{ij} = 0,\ (i,j) \notin \mathcal{E^F}}
\end{mini*}
We propose Algorithm \ref{algo:opt} that solves the above optimization problem \eqref{eqn:P4} in a distributed manner using ADMM for directed graphs.  The steps involved in the derivation are given in the Appendix.\qed

\begin{algorithm} 
  \caption {Local estimation of optimal weight matrix $A_1^{*}$ by agent $i \in \mathcal{F}$ over a directed network using ADMM}
  \label{algo:opt}
  \begin{algorithmic}
  \State \textbf{Initialize :} $\rho > 0$, $(A_1)_i(0) = 0_{n-m}$, $Z_i(0) = 0_{n-m}$, $C_i(0) = 0_{n-m}$, $v_i(0,0)>0$ to a small value$,\ \forall\ i \in \mathcal{F}$ 
    \State \textbf{Iterate :} $\forall\ i \in \mathcal{F}$
    \For{ $k \geq 0$ }
    \State \textbf{Primal Update :}
    \State Evaluate $(A_1)_i(k+1)$ as per \eqref{eqn:Ai}.
    \State Initialize $M_i(0,k) = (A_1)_i(k+1)$
    \For{ $t = 0,1,\hdots,H-1$ }
    \State Send $v_i(t,k)$ and $M_i(t,k)$ to 
    \State out-neighbours.
    \State Compute $v_i(t+1,k)$ using \eqref{eqn:v}
    \State Compute $M_i(t+1,k)$ using \eqref{eqn:M}
    \EndFor
    \State Fix $v_i(0,k+1) = v_i(H,k)$
    \State Fix $Z_i(k+1) = M_i(H,k)$
    \State \textbf{Dual Update :}
    \State Compute $C_i(k+1)$ as per \eqref{eqn:C}.
    \EndFor  
  \end{algorithmic}
\end{algorithm}
\begin{remark}\label{rem:avgestimate}
The update step for the primal variable $Z$ in optimization problem \eqref{eqn:P4} involves computation of an average across all agents over a directed network, as shown in \eqref{eqn:Zi}. Hence for each iteration $k$ in Algorithm 2, every agent $i \in \mathcal{F}$ computes an estimate of $Z$, denoted as $Z_i$, by performing $H$ rounds of communication in the inner loop. The performance of the algorithm improves as $H$ increases. A lower bound on $H$ for the convergence of the algorithm is discussed in \cite{rokad2020distributed}.
\end{remark}

\begin{remark} ({\em Stopping Criterion for Algorithm \ref{algo:opt}}) Each of the agent $i \in \mathcal{F}$ fixes a small, arbitrary tolerance value, $\epsilon > 0$, and calculates its \textit{residual} values at each iteration $k$ as
\begin{equation}
    r_{ij}^1(k) = \frac{1}{n-m}\big\lVert (A_1)_i(k+1)-(A_1)_j(k+1) \big\rVert_F,\ j \in \mathcal{N}^{in}_i
\end{equation}
\begin{equation}
    r_{ij}^2(k) = \big|(A_1)_i^{ij}(k)\big|,\ j \notin \mathcal{N}^{in}_i
\end{equation}
which indicates the extent of violations of the constraints involved. The stopping criterion for Algorithm \ref{algo:opt} is when
\begin{equation}
    R_i(k) = max\big\{r_{ij}^1(k),\ r_{ij}^2(k)\ ;\forall\ j \big\} \leq \epsilon\ \forall\ i \in \mathcal{F}
\end{equation}
\end{remark}

\section{EXAMPLE}
In this section, we present a numerical example that illustrates our result when all the agents are in 2-D plane.

\begin{figure}
  \centerline{\includegraphics[scale=0.5]{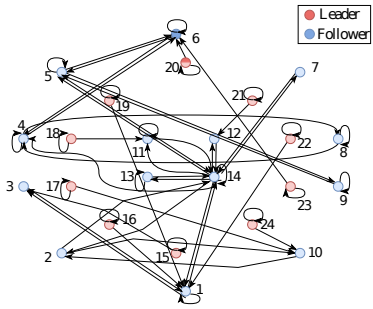}}
   \caption{A directed topology of $n=24$ agents}
   \label{fig:grp}
\end{figure}

\begin{figure}
  \centerline{\includegraphics[scale=0.35]{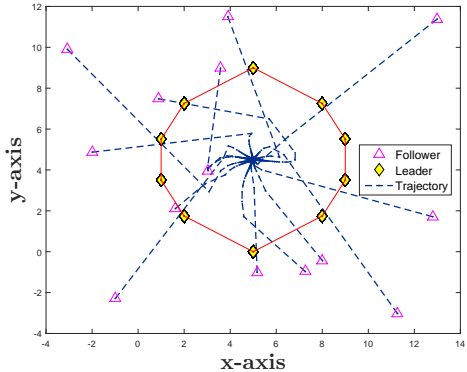}}
  \caption{State trajectories of all agents}
  \label{fig:pplotFig}
\end{figure}

We consider a set of 24 agents interacting over a directed network topology, i.e., $n=24$. Assume that there are $m=10$ leaders, and the rest are followers interacting as shown in  Fig. \ref{fig:grp}.
Let $x^1$ and $x^2$ denote two sets of state variable that indicates the position of each agent in 2-D plane.
The  state values of the stationary leader agents are  given by  $x^1_L = $\big[5 3 2 2 3 5 7 8 8 7\big]$^T$ and $x^2_L =$\big[1 2 4 5 7 8 7 5 4 2\big]$^T$. 
Algorithm \ref{algo:opt} is employed to compute an optimal $A_1$ matrix in a distributed manner.  Then using that $A_1$ and any random choice of $A_2$ following Assumption 1, 
Algorithm \ref{algo:avgcont} is employed to ensure that follower agents reach a consensus to the average of the stationary leader states which is $(\bar{x}^1_L,\bar{x}^2_L) = (5,4.5)$, is illustrated in Fig. \ref{fig:pplotFig}.

\begin{table}[ht!] 
  \begin{center}
    \begin{tabular}{c|c|c}
      \hline 
      Computation method & $\lVert A_1 - \frac{1}{n-m}\mathbf{11}^T \rVert_2$ & $\rho(A_1 - \frac{1}{n-m}\mathbf{11}^T)$ \\
      \hline
      Centralized $\lbrack A_1^{*} \rbrack$ &  \textbf{0.7071} & 0.5010\\
      WBA & 0.9344 & 0.9317\\
      Algorithm \ref{algo:opt} &  \textbf{0.7086} & 0.5146\\
      \hline
    \end{tabular}
    \caption{Convergence factors of different weight matrices for the network in Fig.\ref{fig:grp}}
    \label{tab:table}
  \end{center}
  \vspace{-0.6cm}
\end{table}
\begin{figure}[ht!]
  \centerline{\includegraphics[scale=0.35]{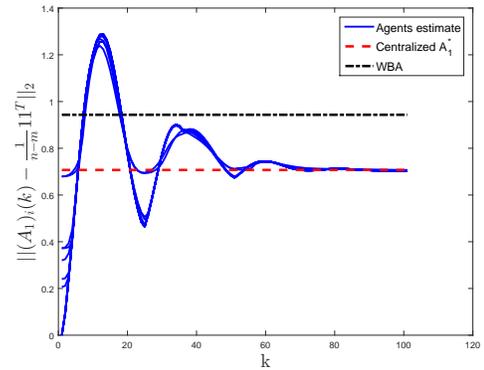}}
  \caption{Objective variation plot of the agents to the centralized optimal value $\lVert A_1^{*} - \frac{1}{n-m}\mathbf{11}^T \rVert_2$}
  \label{fig:ovplotFig}
\end{figure}
The convergence factors of the weight matrices $(A_1)$ generated by different methods have been given in Table \ref{tab:table}, for comparison purposes. A weight-balance algorithm (WBA) proposed in \cite{makhdoumi2015graph} generates balancing weights for a directed network, which, although not optimal, is computed in a distributed manner.
\begin{figure}[ht!]
  \centerline{\includegraphics[scale=0.35]{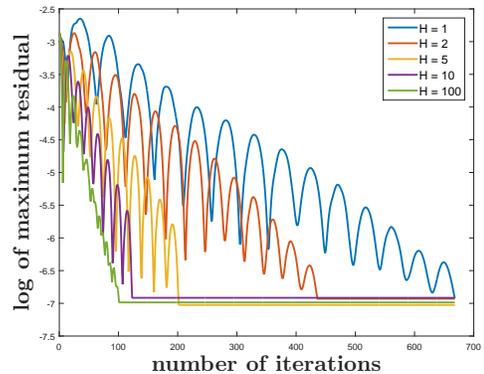}}
  \caption{Maximum log residual plot across all agents with $\rho = 5$}
  \label{fig:rplotFig}
\end{figure}
\\It can be observed that $A_1$ computed by Algorithm \ref{algo:opt} is quite close to the optimal $A_1^{*}$ given by centralized methods as depicted in Fig. \ref{fig:ovplotFig}, and performs much better than that of WBA. In Remark \ref{rem:avgestimate}, we explained how the performance of Algorithm 2 improves with increase in $H$.  This is illustrated in Fig. \ref{fig:rplotFig} where we observe that the \textit{residual} values decreases faster for increase in values of $H$.

\section{CONCLUSION}
We considered a special case of the discrete-time containment control problem for multi-agent systems with multiple stationary leaders, interacting over a directed communication topology. 
We studied the problem of achieving faster convergence to a consensus among followers to the average of the leaders' state, terming it as the average consensus containment. First, we provided sufficient conditions that guarantee average consensus containment. 
A distributed algorithm 
was proposed that ensured average consensus containment. Using this analysis, an optimization problem was formulated to maximize the convergence rate which was also solved in a  distributed manner using ADMM algorithm.
Through numerical examples, we  demonstrated that the optimal convergence rate obtained using locally computed weights in a distributed manner is very close to the optimal value obtained by centralized methods.    
Moving forward, it would be interesting to investigate analogous results using dynamic leaders, and under a continuous-time setting. 
\vspace{-0.13cm}

\bibliographystyle{plain}
\bibliography{Bibliography.bib}
\section*{Appendix}
\section*{Proof of Theorem 1}
Since $A_1$ is a column stochastic matrix associated with a strongly-connected graph (Assumption 1), it is irreducible \cite{ding2011equality}. 
Furthermore, \cite[Sec. 8.3]{meyer2000matrix},  $A_1$ is primitive as its non-negative irreducible with positive diagonal entries due to each follower knowing its own state. Hence, $1$ is a simple eigenvalue of $A_1$. Let $u_r$ and $u_l$ respectively denote the right and left eigenvector of $A_1$ corresponding to eigenvalue  $1$. Note that $u_l=\frac{1_{n-m}}{\sqrt{n-m}}$.
Using \cite[Sec 7.10]{meyer2000matrix} we get
\begin{equation} \underset{k \rightarrow \infty}{\lim}A_1^k = G, \null \label{eqn:lim} 
\end{equation} 
where, $G$ is the spectral projector onto nullspace of $(I_{n-m}-A_1)$ along the columnspace of $(I_{n-m}-A_1)$, and is given by $G=u_ru_l^T$.
Also, as per \cite[Sec 8.4]{meyer2000matrix},  $A_1$ being column-stochastic, is Ces\'{a}ro summable, with the Ces\'{a}ro limit \\
\begin{equation}
\underset{k \rightarrow \infty}{\lim} \frac{A_1^{k-1}+...+A_1+I_{n-m}}{k} = G, \null \label{eqn:cesaro}\\\\
\end{equation}
From \eqref{eqn:lim}, \eqref{eqn:cesaro},\\ 
\begin{equation}
\underset{k \rightarrow \infty}{\lim} \frac{A_1^{k-1}+...+A_1+I_{n-m}}{k} = u_ru_l^T \null \label{eqn:limces} \\\\
\end{equation}
In Algorithm \ref{algo:avgcont}, from Step 3 we get
$x_i(k) = \frac{\left[s(k)\right]_i}{\left[w(k)\right]_i}$ which when using the update law $\forall\ i \in \mathcal{F}$, becomes\\ $x_i(k) = \frac{\left[A_1^kx_F(0) + \{A_1^{k-1}+A_1^{k-2}+...+A_1+I_{n-m}\}A_2x_L(0)\right]_i}{\left[A_1^k 1_{n-m} + \{A_1^{k-1}+A_1^{k-2}+...+A_1+I_{n-m}\}A_2 1_m\right]_i}$\\\\
Dividing  numerator and denominator by $k$ and taking limits  with $x_i^{s} := \underset{k \rightarrow \infty}{\lim}x_i(k)$.\\
$x_i^{s}= \frac{\left[\underset{k \rightarrow \infty}{\lim} \frac{A_1^kx_F(0)}{k} + \underset{k \rightarrow \infty}{\lim} \big\{\frac{A_1^{k-1}+A_1^{k-2}+...+A_1+I_{n-m}}{k}\big\}A_2x_L(0)\right]_i}{\left[\underset{k \rightarrow \infty}{\lim} \frac{A_1^k 1_{n-m}}{k} + \underset{k \rightarrow \infty}{\lim} \big\{\frac{A_1^{k-1}+A_1^{k-2}+...+A_1+I_{n-m}}{k}\big\}A_2 1_m\right]_i}$\\
From \eqref{eqn:lim},  $\left[\underset{k \rightarrow \infty}{\lim} \frac{A_1^kx_F(0)}{k}\right]_i \rightarrow 0$ and $\left[\underset{k \rightarrow \infty}{\lim} \frac{A_1^k 1_{n-m}}{k}\right]_i \rightarrow 0$.
\\From \eqref{eqn:lim}, \eqref{eqn:limces}, we get 
\begin{equation}
    x_i^{s} = \frac{\left[\underset{k \rightarrow \infty}{\lim}A_1^kA_2x_L(0) \right]_i}{\left[\underset{k \rightarrow \infty}{\lim}A_1^kA_2 1_m \right]_i} = \frac{\left[u_r 1_{n-m}^TA_2x_L(0) \right]_i}{\left[u_r 1_{n-m}^TA_2 1_m \right]_i}
    \label{eqn:xs}
\end{equation}
Since, $A_2$ is column stochastic, $1_{n-m}^TA_2 = 1_m^{T}$. Hence,
\[x_i^{s} = \frac{\left[u_r 1_m^{T} x_L(0) \right]_i}{\left[u_r 1_m^{T} 1_m \right]_i} = \frac{\left[u_r 1_m^{T} x_L(0) \right]_i}{\left[ u_rm\right]_i}= \frac{1_m^{T} x_L(0)}{m} = :\overline{x}_L(0)\] 
From the above equation, we observe that the followers are in consensus at the average of the leaders' state values. Hence the proof.\qed

\section*{Proof of Lemma 1}
From \eqref{eqn:xs} we have $\quad \forall\ i \in \mathcal{F}$, 
\begin{equation}
   \underset{k \rightarrow \infty}{\lim}x_i(k) = \frac{\left[\underset{k \rightarrow \infty}{\lim}A_1^kA_2x_L(0) \right]_i}{\left[\underset{k \rightarrow \infty}{\lim}A_1^kA_2 1_m \right]_i}
    \label{eqn:xis}
\end{equation}
$A_1$ is primitive and as such, its eigenvalues satisfy $\lambda_1=1 > \lambda_2 \geq \hdots \geq \lambda_{n-m}$. Therefore, the convergence rate of $x_i(k)$ is determined by $\lambda_2 = \rho\big(A_1 - \frac{1}{n-m}\mathbf{11}^T\big)$.\qed

\section*{Derivation of Algorithm \ref{algo:opt}}
The augmented Lagrangian ($L_{\rho}$) for \eqref{eqn:P4} is given by 
\begin{equation}
    \begin{aligned}
        L_{\rho} = &\ \frac{1}{n-m}\overset{n-m}{\underset{i=1}{\Sigma}}\bigg\lVert (A_1)_i-\frac{1}{n-m}\mathbf{11}^T \bigg\rVert_2\  +\\
        &\ \overset{n-m}{\underset{i=1}{\Sigma}} trace\big[ ((A_1)_i - Z)^T C_i\big] + \overset{n-m}{\underset{i=1}{\Sigma}}\frac{\rho}{2} \bigg\lVert (A_1)_i - Z \bigg\rVert_F^2
    \end{aligned}
    \label{eqn:Lp}
\end{equation}
where $\rho > 0$ is the penalty parameter. Now, the standard 2-block ADMM algorithm for \eqref{eqn:P4} is given below. 
The primal variables $\big[(A_1)_i(k)\ \forall\ i \in \mathcal{F},\ Z(k)\big]$ are updated sequentially and dual variables $\big[C_i(k)\ \forall\ i \in \mathcal{F}\big]$ are updated afterwards. 
The constraint set $S:= \big\{(A_1)_i\ |\ (A_1)_i\mathbf{1}=\mathbf{1},\ (A_1)_i^T\mathbf{1}=\mathbf{1},\ (A_1)_i^{ij}=0,\ (i,j) \notin \mathcal{E^\mathcal{F}}\big\}$ has been incorporated into primal update steps rather than being dualized.\\
\textbf{$Z$ update step :} $Z(k+1) = \underset{Z}{argmin} \big\{L_{\rho}\big\}$
\begin{equation}
    \begin{aligned}
        Z(k+1) =&\ \frac{1}{n-m}\overset{n-m}{\underset{i=1}{\Sigma}}\bigg[ \frac{C_i(k)}{\rho} + (A_1)_i(k+1) \bigg]
    \end{aligned}
    \label{eqn:Zi}
\end{equation}
The $Z$ update step involves computing an average across all follower agents over a directed network. Hence we use an approach based on \textit{dynamic average consensus}, proposed in \cite{rokad2020distributed}, to get the $Z$ update step, in a distributed manner. Each agent $i \in \mathcal{F}$ maintains $Z_i$ which is an estimate for $Z$. Before we provide a distributed update law for $Z$, we write the update law for the other variables. \\
\textbf{$(A_1)_i$ update step :} $(A_1)_i(k+1) = \underset{S}{argmin} \{L_{\rho}\}$
\begin{equation}
    \begin{aligned}
        (&A_1)_i(k+1) = \underset{S}{argmin}\Bigg\{\frac{\lVert(A_1)_i-\frac{1}{n-m}\mathbf{11}^T\rVert_2}{n-m}\ +\\ &\ trace\big[ ((A_1)_i - Z_i(k))^T C_i(k)\big]\ + \frac{\rho}{2} \bigg\lVert (A_1)_i - Z_i(k) \bigg\rVert_F^2\Bigg\}
    \end{aligned}
    \label{eqn:Ai}
\end{equation}
\textbf{$C_i$ update step :}
\begin{equation}
C_i(k+1) = C_i(k) + \rho \big[ (A_1)_i(k+1) - Z_i(k+1) \big] \label{eqn:C}    \end{equation}
Note that in both \eqref{eqn:Ai} and \eqref{eqn:C}, we use the estimate $Z_i$ instead of $Z$. Next, we simplify the $Z$ update step in \eqref{eqn:Zi}.
At every iteration $k$, each agent initializes a variable $M_i(0,k): = Z_i(k) + (A_1)_i(k+1) - (A_1)_i(k) + \frac{C_i(k) - C_i(k-1)}{\rho}$. From \eqref{eqn:C}, this is further simplified as $M_i(0,k) = (A_1)_i(k+1)$. For some $H\geq 1$, $0\leq t\leq H-1$, $i \in \mathcal{F}$ do\\
\begin{equation}
    \begin{aligned}
        M_i(t+1,k) =&\ \big[ 1-d_iv_i(t,k) \big] M_i(t,k)\ +\\ &\ \underset{j \in \mathcal{N}^{in}_i}{\Sigma}v_j(t,k)M_j(t,k)
    \end{aligned}
    \label{eqn:M}
\end{equation}
\begin{equation}
    v_i(t+1,k) = \frac{1}{2} \big[ v_i(t,k) + \frac{1}{d_i} \underset{j \in \mathcal{N}^{in}_i}{\Sigma} v_j(t,k) \big]
    \label{eqn:v}
\end{equation}
Then $Z_i(k+1)=M_i(H,k)$. In \eqref{eqn:M}, $d_i := |\mathcal{N}^{out}_i|$ is the out-degree of $i^{th}$ follower agent. In \eqref{eqn:v}, $v_i(t,k) \in \mathbb{R}$ is the node-weight used by $i^{th}$ agent to scale its outgoing information and evolves dynamically as shown. 
As defined in \cite{makhdoumi2015graph}, the node-weights are initialized as $v_i(0,0) \leq (1/d^{*})^{2D+1}$, where $D$ is the diameter of the subgraph $\mathcal{G^F}$ and $d^{*} = \underset{i \in \mathcal{F}}{max}\ d_i$. Each of the agent $i \in \mathcal{F}$ initialises a small value as its node-weight, $v_i(0,0) > 0$.\qed

\balance

\end{document}